\documentclass[12pt,a4paper,reqno]{amsart}
\usepackage{amssymb}
\usepackage{amsfonts}
\usepackage{amsthm}
\usepackage{amsmath}
\usepackage{mathrsfs}
\usepackage{amscd}
\usepackage{calc}
\usepackage[latin2]{inputenc}
\usepackage{t1enc}
\usepackage[dvipsnames]{xcolor}
\usepackage{indentfirst}
\usepackage{graphicx}
\usepackage{graphics}
\usepackage{pict2e}
\usepackage{epic}
\numberwithin{equation}{section}
\usepackage[margin=2.9cm]{geometry}
\usepackage{epstopdf} 
\usepackage{mathrsfs}
\usepackage[shortlabels]{enumitem}
\usepackage{amsrefs}

\theoremstyle{plain}
\newtheorem{Th}{Theorem}[section]

\newtheorem{Prop}[Th]{Proposition}

\theoremstyle{definition}
\newtheorem{Def}[Th]{Definition}

\newcommand{\R}{\mathbb{R}}

\renewcommand{\epsilon}{\varepsilon}

\newcommand{\abs}[1]{\lvert #1 \rvert}

\DeclareMathOperator{\supp}{supp}

\title{Scaling Asymptotics for Ladder Sequences of Spherical Harmonics at Caustic Latitudes}
\author[M. Geis]{Michael Geis}
\address{Department of Mathematics, Northwestern University, 2033 Sheridan Road, Evanston, IL 60208}
\email{mlg@math.northwestern.edu}
\thanks{Written in progress of a PhD in mathematics at Northwestern University}
\date{August 3, 2022}

\begin{document}

\maketitle

\begin{abstract}
We study the concentration of ladder sequences of spherical harmonics on caustic latitude circles. We prove that they have Airy scaling asymptotics. We also determine the weak* limit of certain empirical measures of $L^2$ norms of restrictions of spherical harmonics to these latitude circles.
\end{abstract}

\section{Introduction}

Consider the round sphere $(S^2,g_{can})$ with standard polar coordinates $(\phi,\theta) \in (0,\pi) \times (0,2\pi)$ where $\theta$ is the polar angle measured relative to fixed meridian geodesic. We study scaling asymptotics of certain sequences of the standard $L^2$ normalized spherical harmonics,

\begin{equation}
    Y^m_N(\phi,\theta) = \sqrt{\frac{2N + 1}{4\pi} \frac{(N - m)!}{(N + m)!}} P^m_N(\cos\phi)e^{im\theta},
\end{equation}

in which $m_k,N_k \to \infty$ while the ratio $c = m_k/(N_k +\frac{1}{2})$ is held fixed. These are called rational ladder sequences and have the special feature that they are semi-classical Lagrangian distributions which concentrate on Lagrangian tori $T_c \subset S^*S^2$. The torus $T_c$ projects to an annular band around the equator $\gamma_e = \{\phi = \frac{\pi}{2}\}$ and the projection has a fold singularity over the bounding latitude circles, $\gamma^\pm_c$ determined by $\sin \phi = c$. We obtain scaling asymptotics in an $(N + \frac{1}{2})^{-\frac{2}{3}}$ neighborhood around the caustic latitude circles

\begin{Th}
There exists an $\epsilon > 0$ such that if $x = (\phi,\theta)$ with $c < \sin\phi < c + \epsilon$ then, letting $h_k = \left(N_k + \frac{1}{2}\right)^{-1}$, there are a sequence of smooth half densities $u_{i,j}$ on $S^2$ such that 

\begin{equation}
\begin{aligned}
    Y^{m_k}_{N_k}(x) \sqrt{dV_g(x)} \sim& Ai\left(-h_k^{-\frac{2}{3}} \rho(x) \right) \sum_{n = 0}^\infty u_{0,n}(x) h_k^{-{\frac{1}{6}} + n} +  \\ & Ai'\left(h_k^{-\frac{2}{3}} \rho(x)\right)\sum_{n = 0} u_{1,n}(x) h_k^{\frac{1}{6}+n} 
    \end{aligned}
\end{equation}

The argument of the Airy function and its derivative is 

\begin{equation}
    \rho(x) = \left( \frac{4}{3} \int_{\gamma_x} \alpha \right)^\frac{2}{3}
\end{equation}

Here, $\gamma_x$ is the geodesic arc joining the two pre-images $\pi^{-1}(x) \in T_c$ and $\alpha$ is the canonical 1-form on $T^*S^2$. The arc is oriented so as to make the integral positive. The leading order amplitude is given by 

\begin{equation} \label{HALFDENSITYINTERP}
    u_{0,0}(x) =  (2\pi) \rho(x)^\frac{1}{4} \pi_*\sqrt{d\mu_c}
\end{equation}

where $d\mu_c$ for the normalized joint flow invariant density on $T_c$.

\end{Th}

The articles recent work of Galkowski and Toth \cite{GT1},\cite{T1} obtain sharp decay estimates for joint eigenfunctions in the forbidden region $S^2 \setminus \pi(T_c)$. The scaling asymptotics we obtain are only valid on the `allowed region side' of the caustic latitude, however it should be possible to extend them to two sided Airy asymptotics which agree with the decay proven by Galkowski-Toth. It should be possible to extend the Airy scaling results to joint eigenfunctions on a convex surface of revolution. One can still separate variables and conjugate the 1D Sturm-Liouville problem to a Schr\"odinger operator. Once a quasi-mode is constructed, the rest of the argument works the same way. Scaling asymptotics of the Legendre functions $P^m_N$ along such sequences were studied previously in \cite{Th1} \cite{Ol1}, however the ODE methods they applied did not provide explicit expressions for the quantities appearing in the expansion. We re-derive the expansion of the Legendre functions by constructing explicit quasi-modes for the Legendre operator using the well-known Maslov-WKB quantization procedure (\cite{CdV2},\cite{ES},\cite{Z1}) which approximate the Legendre functions locally uniformly up to $O(h^\infty)$ error. The quasi-mode is expressible as an oscillatory integral with a degenerate critical point in the phase near the turning points. From this, the Airy expansion is obtained by putting the phase function in a cubic normal form. This idea was first developed by Chester Friedman and Ursell \cite{CFU} and later by Ludwig \cite{Lud1},\cite{Lud2}. It is also discussed by Guillemin and Sternberg \cite{LPDO1} and H\"ormander \cite{GS1} The advantage of this approach is that it expresses the quantities appearing in the expansion of the Legendre functions in terms of the geometry on $S^2$.

\vspace{5mm}

In the previous paper \cite{MG}, we studied the empirical measures

\begin{equation}
    \mu_N(t) = \sum_{m = -\ell}^\ell ||\varphi^m_\ell||^2_{L^2(H)}\delta_0\left(t - \frac{m}{\ell}\right)
    \end{equation}
    
On a convex surface of revolution $(S^2,g)$, where $H$ is the unique rotationally invariant closed geodesic. It was shown that these measures tend to a weak limit which exhibits $(1 - t^2)^{-\frac{1}{2}}$ type blow-up at the end points due to the fact that the end points $t = \pm 1$ correspond to the Gaussian beam ladder sequences of joint eigenfunctions $(m/\ell = \pm 1)$ which peak on the geodesic $H$.

\begin{Th}
Let $\gamma_{c_0}$ be either of the two latitude circles determined by $\sin \phi = c_0$. The normalized empirical measures,

\begin{equation} \label{EMPIRICALMEASURE}
    \mu_{N,c_0} = \frac{1}{c_0(N + \frac{1}{2})}\sum_{m = -N}^N ||Y^m_N||^2_{L^2(\gamma_{c_0})} \delta_0\left(t - \frac{m}{N}\right) 
\end{equation}

Then for all $f \in C^0([-1,1])$,

\begin{equation}
    \int f(t) \, d\mu_{N,c_0}(t) = \frac{1}{c_0 \pi}\int_{-c_0}^{c_0} f(t) \left(1 - \left(\frac{t_0}{c_0}\right)^2\right)^{-\frac{1}{2}} \, dt
\end{equation}

\end{Th}

This reflects the fact that the latitude $\gamma_{c_0}$ is the caustic curve for the ladder sequence with ratio $m/N \to c$ and these spherical harmonics are largest there, with an $(N + \frac{1}{2})^{1/6}$ peak due to the Airy bump. The measures are supported on $[-c_0,c_0]$ reflecting the fact that when $c = m/N > c_0$, $\gamma_c$ is outside of the projection of the Lagrangian $T_c$ which makes the corresponding ladder sequence $O(h^\infty)$ on $\gamma_{c}$

\subsection*{Acknowledgements}
I would like to thank my advisor, Steve Zelditch for his unwavering support, patience and guidance.

\subsection{Outline}
In section 2, we conjugate the associated Legendre operator to a Schr\"odinger operator on $I = (0,\pi)$ and construct a global WKB quasi-mode (approximate eigenfunction) in such a way that it is a locally uniform approximation to $P^m_N(\cos \phi)$, $\phi \in (0,\pi)$ as $m,N \to \infty$ with the ratio $c = m/(N+\frac{1}{2})$ fixed. Section 3 contains the derivation of the Airy expansion of the quasi-mode. In section 4 we explain the relevant geometry on $S^2$ and connect it to the expansion of the Legendre functions. Section five contains the short calculation of the weak limit of the empirical measures $\eqref{EMPIRICALMEASURE}$. The rest of this section contains background on the Legendre functions and an overview of the relevant machinery of semi-classical Lagrangian distributions for those unfamiliar.

\subsection{Background on Legendre functions}

To establish notation and collect basic facts we quote the following classical results about the Legendre functions and refer to the standard references \cite{Ol1},\cite{Hob} for more detail. We note that these functions are called `Ferrer's functions' or `Legendre functions on the cut' by some authors. For each pair of integers $(m,N)$ with $0 \leq m \leq N$, let $P^m_N(x)$ be the following function defined for $x \in [-1,1]$:

\begin{equation}
    P^m_N(x) := \left(\left(N + \frac{1}{2}\right) \frac{(N - m)!}{(N+m)!}\right)^\frac{1}{2}\frac{1}{2^N N!} (1-x^2)^{\frac{m}{2}} \partial_x^{N+m}(x^2-1)^N
\end{equation}

We refer to $P^m_N(x)$ as the normalized Legendre function of degree $N$ and order $m$. They are real-valued, smooth on $(-1,1)$, and satisfy

\begin{equation} 
    (1-x^2)\partial_x^2P^m_N(x) - 2x\partial_xP^m_N(x) + \left( \left(N + \frac{1}{2}\right)^2 - \frac{m^2}{1-x^2} - \frac{1}{4} \right) P^m_N(x) = 0,
\end{equation}

\begin{equation}\int_{-1}^1 P_N^m(x)^2 \, dx = 1.
\end{equation}

\begin{Prop} \label{LEGSPEC}
For $m \in \mathbb{Z}_{\geq 0}$, define the (positive) Legendre operator $L_m$, 

\begin{equation} \label{LEGOP}
    L_m := -\partial_x(1-x^2)\partial_x + \frac{m^2}{1-x^2} + \frac{1}{4}
\end{equation}

As an unbounded operator on $L^2[-1,1]$ with domain $C_c^\infty([-1,1],dx)$, $L_m$ has only discrete spectrum consisting of simple eigenvalues 

$$
    \text{Spec}(L_m) = \left\{ \left(N + \frac{1}{2}\right)^2 ~|~ N \in \mathbb{N}, N \geq m \right\}.$$
    
    Each eigenspace is the complex span of $P^m_N(x)$ and the set $\{ P_N^m \}_{N = m}^\infty$ is an orthonormal basis of $L^2([-1,1], dx)$

\end{Prop}

For a proof, see \cite{OP}. The formula 

$$P_{n+1}'(x) = xP'_n(x) + (n+1)P_n(x)$$

together with $P^0_N(1) = 1$ implies that for all $0 \leq m \leq N$, $P^m_N(x)$ is positive near $x = 1$. Depending on the relative of parity of $m$ and $N$, $P^m_N(x)$ is either odd or even, $P^m_N(-x) = (-1)^{m + N} P^m_N(x)$. We will use these properties to match the quasi-mode with $P^m_N$.

\subsection{Background on oscillatory functions associated to Lagrangian manifolds}

This section contains a review of the basic theory of oscillatory integrals which we will use in the construction of the quasi-mode in section 2.

\vspace{5mm}

Let $(M^n,g)$ be a Riemannian manifold. The theory reviewed here depends upon working with smooth half densities rather than functions. Fix a smooth, positive density $\nu$ on $M$. We may then identify functions with half densities via the isomorphism

$$f(x) \cong f(x) \sqrt{\nu}$$

Let $\Lambda \subset T^*M$ be a compact Lagrangian submanifold. In order to define the space $\mathscr{O}^*(M,\Lambda)$ of oscillatory half densities associated to $\Lambda$, we fix a locally finite open cover $\{U_j\}$ of $\Lambda$ such that for each $j$, there exists a phase function $\psi_j(x,\theta) \in C^\infty(V_j \times \R^{N_j}, \R)$ defined on some open subsets $V_j \subset M$ which are small enough so that the maps

$$i_{\psi_j} : (x,\theta) \ni C_{\psi_j} \mapsto (x, d_x\psi_j(x,\theta))$$

are embeddings onto $U_j \subset \Lambda$. Here $C_{\psi_j}$ is the zero set of $d_\theta\psi_j$ which is assumed to be an $n$ dimensional submanifold. We further fix a partition of unity $\chi_j$ subordinate to this cover. 

\begin{Def}
The space $\mathscr{O}^\mu(M,\Lambda)$ is the space of all half densities which can be written in the form 

\begin{equation} \label{LOCALOSC}
u(x,h) = \left(\sum_{j} (2\pi h)^{-\frac{N_j}{2}}\int_{\R^{N_j}} a_j(x,\theta,h) e^{\frac{i}{h}\psi_j(x,\theta)} \, d\theta \right)  \sqrt{\nu}
\end{equation}

\begin{equation}
a_j(x,\theta) \sim \sum_{n = 0}^\infty a_{j,n}(x,\theta) h^{\mu + n}    
\end{equation}

\end{Def}

where each $a_j(x,\theta)$ is a smooth function with compact support. We write $\mathscr{O}^{\infty}(M,\Lambda) = \cap_{\mu \in \R} \mathscr{O}^\mu(M,\Lambda)$ and when $h$ is restricted to take values in a particular sequence $h_k$ we will signify this with the notation $\mathscr{O}^\mu(M,\Lambda,h_k)$. Associated to each $u(x,h) \in \mathscr{O}^\mu(M,\Lambda)$ is a geometric object $\sigma(u)$ called its principal symbol, which is a section of a certain line bundle over $\Lambda$. To define it, we first recall that the Maslov bundle $\mathbb{L} \to \Lambda$ is a flat complex line bundle which can be described concretely using the choice of $\{U_j,\psi_j\}$. On $U_i \cap U_j$, define the locally constant functions

$$m_{ij}(\lambda) = \frac{1}{2}\left(\text{Sgn}\,\partial_\theta^2\psi_j - \text{Sgn}\,\partial_\theta^2\psi_i\right)$$

where $\partial_\theta^2\psi$ is the hessian with respect to the fiber variables. The functions $\exp i\frac{\pi}{2}m_{ij}(\lambda)$ are the transition functions of the Maslov bundle on $U_i \cap U_j$. The choice of phase functions determines a canonical section, $s$, of $\mathbb{L}$ by

$$s_j(\lambda) = \exp i\frac{\pi}{4}\text{Sgn} \,  d_\theta^2\psi_j(\lambda) \quad \quad \lambda \in U_j$$

Let $\Psi_j$ be the lift of $\psi_j$ to $U_j$ via the map $i_{\psi_j}$ and $\Omega^{\frac{1}{2}} \to \Lambda$ be the half density bundle over $\Lambda$. Fix a smooth positive density $\rho_0$ on $\Lambda$ and define the space of symbols of order $\mu$, $S^\mu(\Lambda)$, to be the set of all smooth sections of $\Omega^\frac{1}{2}\otimes \mathbb{L} \to \Lambda$ which may be written in the form

\begin{equation} \label{PRINCIPALSYMBOL}
    h^\mu\left( \sum_j \exp i\frac{\Psi_j(\lambda)}{h} f_j(\lambda) s_j(\lambda) + O(h)\right) \sqrt{\rho_0}
\end{equation}

where $f_j$ are smooth functions on $\Lambda$ with $\supp f_j(\lambda) \subset U_j$. The principal symbol map $\sigma : \mathscr{O}^\mu(M,\Lambda) \to S^\mu(\Lambda)/S^{\mu + 1}(\Lambda)$ is defined so that when $u(x,h)$ is written in the form \eqref{LOCALOSC} then

\begin{equation}
    [\sigma(u)](\lambda) = h^{\mu}\left(\sum_j \exp i{\frac{\Psi_j(\lambda)}{h}} a_{j,0}(\lambda) g_j(\lambda)s_j(\lambda)\right) \sqrt{\rho_0}
\end{equation}

Here, the $g_j$ are smooth functions on $U_j$ defined by 

$$g_j \sqrt{\rho_0} = (i_{\psi_j}^{-1})^*\sqrt{d_{C_{\psi_j}}} \quad \quad \lambda \in U_j$$

where $d_{C_{\psi_j}}$ is the canonical $\delta$-density on the critical set $C_{\psi_j}$ determined by the density $\nu \otimes \abs{d\theta}^\frac{1}{2}$ on $V_j \times \R^{N_j}$. Next, we a map which takes a symbol to an oscillatory half density, 

$$\mathscr{Q} : S^\mu(\Lambda) \to \mathscr{O}^\mu(M,\Lambda).$$

Suppose that $\sigma \in S^\mu$ is written in the form \eqref{PRINCIPALSYMBOL} (which is always possible using $\chi_j$). Define $\mathscr{Q}(\sigma)$ to be the smooth half density \eqref{LOCALOSC} with amplitudes $a_j$ chosen so that $(i_{\psi_j}^{-1})^{*}a_j g_j = f_j$. Then $\mathscr{Q}$ is a right inverse for the principal symbol map. It depends on the choices $(U_j,\psi_j,\chi_j)$ while the principal symbol map does not. Finally, suppose that $P$ is an order zero semi-classical pseudo-differential operator with principal symbol $p_0$ and with sub-principal symbol equal to zero. If $p_0 = 0$ on $\Lambda$ and $\rho$ is a density on $\Lambda$ invariant under the Hamiltonian flow $t \mapsto \exp tX_{p_0}$ of $p_0$, then for any $u \in \mathscr{O}^\mu(M,\Lambda)$, $Pu \in \mathscr{O}^{\mu+1}(M,\Lambda)$ and if $u(x,h)$ has principal symbol

\begin{equation}
    \sigma(u) = \left( \sum_j \exp i \frac{\Psi_j(\lambda)}{h_k} f_j(\lambda) s_j(\lambda) \right) \sqrt{\rho}
\end{equation}

then the order $\mu+1$ symbol of $Pu$ is

\begin{equation} \label{SUBPRINICPALTERM}
    \sigma(Pu) = \left( \sum_j \exp i\frac{\Psi_j(\lambda)}{h_k} \frac{2}{i}X_{p_0}f_j(\lambda) s_j(\lambda) \right)\sqrt{\rho}
\end{equation}

\section{WKB for the Legendre operator}

We begin by conjugating the Legendre operator on $[-1,1]$ to a Schr\"odinger operator on $I = (0,\pi)$. The following proposition is a straightforward calculation.

\begin{Prop} \label{LEGCONJ}

Let $U$ be the unitary map $U : L^2((-1,1), dx) \to L^2((0,\pi), d\phi)$ 

$$(Uf)(\phi) =  f(\cos \phi) \sqrt{\sin \phi} $$ 


Let $0<c<1$ and define the operator $H_{h,c}$ for $f \in C^\infty((0,\pi))$

\begin{equation}
  H_{h,c}f(\phi) := -h^2f''(\phi) + \left(\frac{c^2}{\sin^2\phi} - \frac{h^2}{4\sin^2\phi}\right)f(\phi)
\end{equation}

Suppose that $m(h)$ is an integer such that $c = m(h)h \in (0,1)$ for all $h$. Then 

$$ h^2UL_{m(h)}U^* = H_{h,c}$$

\end{Prop}





For the remainder of this section we fix once and for all some $c \in (0,1)$ and a rational ladder sequence, that is, integers $0 \leq m_k \leq N_k$ such that for all $k$, $m_k/(N_k + \frac{1}{2}) = c$.

Putting $h_k = (N_k + \frac{1}{2})^{-1}$, it follows from propositions \ref{LEGSPEC} and \ref{LEGCONJ} that the spectrum of $H_{h_k,c}$ is

$$\text{Spec}(H_{h_k,c}) = \left\{ h_k^2\left(N+\frac{1}{2}\right)^2 ~|~ N \geq m_k \right\}.$$

In particular, $1$ is an eigenvalue of $H_{h_k,c}$ for all $k$. Moreover $\ker \left( H_{h_k,c} - 1\right)$ is one dimensional and spanned by $u_{h_k}(\phi) := UP^{m_k}_{N_k}$. It follows that there exists $\delta > 0$ so that $H_{h_k,c}$ has the spectral gap,

\begin{equation} \label{SPECGAP}
    \inf_{\lambda \in \text{Spec}(H_{h_k,c}) \setminus \{ 1 \}} |1 - \lambda| \geq \delta h_k
\end{equation}

\subsection{Construction of a global $h^\infty$ quasi-mode for $H_{h_k,c}$}

We say that a smooth function $v_{h}$ on $I = (0,\pi)$ is a quasi-mode of order $h^\infty$ for $H_{h,c}$ with quasi-eigenvalue $E(h)$ if 

\begin{equation}
    ||\left(H_{h,c} - E(h)\right)v_{h}||_{L^2(I)} = O(h^\infty)
\end{equation}

where $E(h)$ has the semi-classical expansion $E(h) \sim E_0 + \sum_{j=1}^\infty h^jE_j$. Let $(\phi,\tau)$ be coordinates for $T^*\R{}$, $p : T^*\R \to \R$ the natural projection, and

$$f(\phi,\tau) = \tau^2 + \frac{c^2}{\sin^2\phi}$$

be the principal symbol of $H_{h,c}$. The energy curve

\begin{equation}
    \Sigma = \{ f(\phi,\tau) = 1 \}
\end{equation}

is a smooth, closed curve symmetric about $\tau \mapsto -\tau$, intersecting $\{\tau = 0\}$ at $\phi_\pm = \pi/2 \pm \phi_0$ where $\phi_\pm$ are the two solutions of $\sin \phi = c$, $\phi \in (0,\pi)$. We follow the well-known procedure of WKB-Maslov quantization in order to construct a quasi-mode $v_{h_k}$ approximating $u_{h_k} = U P^{m_k}_{N_k}$ locally uniformly on $I$. For the remainder of the section we identify smooth functions on $I$ with smooth half densities on $I$ by 

$$f(\phi) \cong f(\phi)\, \abs{d\phi}^\frac{1}{2}.$$

In this way we may speak of oscillatory functions instead of half densities and we do this without further comment. The rest of this section contains the proof of the following:

\begin{Prop} \label{QMPROP}
There exists a smooth, real-valued function $v_{h_k}(\phi) \in \mathscr{O}^0(I,\Sigma,h_k)$ with $||v_{h_k}||_{L^2(I)} = 1$ and a sequence of real numbers $E_j$ so that if $E(h_k) \sim 1 + h_k^2E_2 + h_k^3E_3 + \cdots$, then

\begin{equation} \label{QMEQN}
||(H_{h_k,c} - E(h_k))v_{h_k}||_{L^2(I)} = O(h_k^\infty).\end{equation}

\vspace{3mm}

Moreover, for any fixed $\phi \in I$, 

\begin{equation} \label{ALLOWEDPTWISE}
    v_{h_k}(\phi) = \begin{cases}
    \sqrt{\frac{2\sin \phi}{\pi}} \ \frac{\cos \left( \frac{1}{h_k}\int_{\gamma_\phi} \alpha + \frac{\pi}{4} \right)}{\left( \sin^2\phi - c^2 \right)^{\frac{1}{4}}} + O(h)_{L^2} \quad N_k - m_k ~\text{odd} \vspace{5mm} \\
    
    -\sqrt{\frac{2\sin \phi}{\pi}} \ \frac{\sin \left( \frac{1}{h_k}\int_{\gamma_\phi} \alpha + \frac{\pi}{4} \right)}{\left( \sin^2\phi - c^2 \right)^{\frac{1}{4}}} + O(h)_{L^2} \quad N_k - m_k ~\text{even}
    \end{cases} 
\end{equation}

and there exists an $\epsilon > 0$ such that if $\phi > \phi_+ - \epsilon$, then 

\begin{equation} \label{TURNINGPTWISE}
    v_{h_k}(\phi) = (2\pi h_k)^{-\frac{1}{2}} \int a(\tau,h) e^{\frac{i}{h_k}(\phi\tau - G_4(\phi))} \, d\tau + O(h^\infty)_{L^2},
\end{equation}

Where $a(\tau,h) \sim \sum_j a_j(\tau) h^j$, $a_0(\tau) = \frac{1}{\sqrt{\pi}} V'(G_4'(\tau))^{-\frac{1}{2}}$, and $G_4(\tau)$ satisfies $G_4(0) = 0$, $(G'_4(\tau),\tau) \in \Sigma$ on the support of $a$.

\end{Prop}

The existence of $O(h^\infty)$ quasi-modes is well known, see for instance \cite{CdV2}, \cite{Du2}, \cite{ES}. The solvability of \eqref{QMEQN} up to error $O(h^\infty)$ requires $(\Sigma, h_k)$ to have the following three properties:  

\begin{Prop} \label{NECESSARYCONDITIONSQUANT}
Let $\alpha = \tau d\phi|_{\Sigma}$, $[\mathfrak{m}] \in H^1(\Sigma, \mathbb{Z})$ be the Maslov class, and $X_f$ be the Hamiltonian vector field of $f$. 

\begin{enumerate}[label = (\alph*)]
    \item For all $k$ large enough, 

\begin{equation} \label{BSCCOHOMOLOGY}
    \frac{1}{2\pi h_k}[\alpha] - \frac{1}{4}[\mathfrak{m}] \in H^1(\Sigma,\mathbb{Z})
\end{equation}

\item There exists a positive density $\rho_0$ invariant under the flow of $X_f$

\item For each smooth function $r_0$ on $\Sigma$ satisfying $\int_{\Sigma} r_0 \, \rho_0 = 0$, there exists a smooth function $r_1$ so that $dr_1(X_f) = r_0$.

\end{enumerate}

\end{Prop}

\begin{proof}

\begin{enumerate}[label = (\alph*)]
    \item Since $\Sigma$ is a curve, we can check this by integration. We define the Maslov class below, but to check this it suffices to know that $\int_{\Sigma}[\mathfrak{m}] = 2$ when $\Sigma$ is oriented counter-clockwise. Since the integral of $\alpha$ is the area enclosed by $\Sigma$ and $\Sigma$ is symmetric across the lines $\phi = \pi/2$, $\tau = 0$, the integral is four times the area of the upper right quadrant, 
    
    $$\int_{\Sigma} \tau \, d\phi = 4c \int_0^{\sqrt{1 - c^2}} \frac{\tau^2 \, d\tau}{(1- \tau^2)\sqrt{1 - c^2 - \tau^2}}  = 2\pi(1 - c)$$

Therefore 

$$\frac{1}{2\pi h_k}\int_{\Sigma}\alpha - \frac{1}{2} = N_k - m_k \in \mathbb{Z}$$


\item The map 

$$i : [0,\pi) \to \Sigma \quad \quad i(t) = \exp tX_f(\phi_+,0)$$ 

is a surjective Lagrangian immersion. To see this, one only needs to note that the period of the Hamiltonian flow through $(\phi_+,0)$ is $\pi$. This follows from the fact the curve $\exp\frac{t}{2}X_f(\phi_+,0)$ can be identified with a geodesic on $S^2$ (See section 4). The density $\rho_0$ defined by $i^*\rho_0 = \pi^{-1}\abs{dt}$ is clearly positive and invariant. 

\item Pulling back under $i$, we may assume $r_0(t)$ is smooth on $[0,\pi)$, $\int_0^\pi r(t) \, \abs{dt} = 0$, and $\lim_{t \to \pi}r_0(t) = r_0(0)$. Then the function $r_1(t) = \int_0^{t} r_0(s) \, \abs{ds}$ solves the equation.

\end{enumerate}

\end{proof}

\subsubsection{Explicit choice of phases and canonical operator}

Let $\{U_j\}_{j=1}^4$ be the open cover of $\Sigma$ described as follows: pick $0 <  s < \frac{\pi}{2}$. and let $U_1 = p^{-1}(\frac{\pi}{2} - s, \frac{\pi}{2} + s) \cap \{ \tau > 0 \}$ and $U_3 = p^{-1}(\frac{\pi}{2} - s, \frac{\pi}{2} + s) \cap \{ \tau < 0 \}$. Then we let $U_2 = \{(\phi,\tau) \in \Sigma ~|~ \phi \in (\phi_-, \frac{\pi}{2} - s + \epsilon \}$ and $U_4 = \{(\phi,\tau) \in \Sigma ~|~ \phi \in (\frac{\pi}{2} + s - \epsilon, \phi_+) \}$ where $0 < \epsilon < s$ can be arbitrary. The sets $U_1,U_3$ are symmetric about $(\phi,\tau) \mapsto (\phi,-\tau)$ and $U_2,U_4$ are symmetric with respect to reflection over $\phi = \frac{\pi}{2}$.

\vspace{10mm}

Let $\chi_j$ be a partition of unity subordinate to this cover so that $\chi_1(\phi,\tau) = \chi_3(\phi,-\tau)$ and $\chi_2(\phi,\tau) = \chi_4(\pi - \phi,\tau)$. We choose local phase functions parametrizing this open cover as follows. For $j = 2,4$ we put

\begin{equation}
    \psi_j(\phi,\tau) = \phi\tau - G_j(\tau),
\end{equation}

where $G_j$ are chosen so that $G_j(0) = 0$ and $(G'_j(\tau),\tau) \in \Sigma$ on the $\tau$-projection of $U_j$. For $\phi \in p(U_1) = p(U_3)$, let

\begin{equation}
    \psi_1(\phi) = \int_{\gamma_\phi} \alpha \quad\quad\quad \psi_3(\phi) = \int_{-\gamma_\phi} \alpha = -\psi_1(\phi)
\end{equation}

Here, $\gamma_\phi$ is the arc joining the turning point $(\phi_+,0)$ to the point $(\phi,\tau) \in U_1$ and $-\gamma_\phi$ is the arc joining $(\phi_+,0)$ to $(\phi,\tau) \in U_3$. Since the lifts $\Psi_j$ of the phases to $\Sigma$ are primitives of $\alpha$, they differ by a constant $\Psi_i - \Psi_j := C_{ij}$ on each $U_i \cap U_j$. It is easy to see for this choice of phases that $C_{12} = C_{23} := C$ and $C_{34} = C_{41} = 0$. Note that this means $\int_{\Sigma} \alpha = 2C$ where the integral is in the counter-clockwise direction. As described in section 1.3, this choice of phases shows that the co-cycle which defines the Maslov class $[\mathfrak{m}]$ is $m_{21} = m_{32} = m_{43} = m_{14} =  \frac{1}{2}$.

\begin{Prop} \label{PHASEPROP}
Define constants $\beta_j$ as follows

$$
\begin{cases}
\beta_1 = -\beta_3 = \frac{\pi}{4} \ \quad N_k - m_k ~\text{odd} \\
\beta_1 = -\beta_3 = \frac{3\pi}{4} \quad N_k - m_k ~\text{even}
\end{cases}
\quad \quad \begin{cases}
\beta_2 = \beta_4 = 0 \quad \quad N_k - m_k ~\text{odd} \\
\beta_2 = 0,\beta_4 = \pi \quad N_k - m_k ~\text{even}
\end{cases}$$

\vspace{2mm}

Then the local expressions

\begin{equation}
    S_j(\lambda) = \exp \, i \left( \frac{\Psi_j(\lambda)}{h_k} + \beta_j  \right)s_{\psi_j}(\lambda) \quad \quad \lambda \in U_j
\end{equation}

define a global section of the Maslov bundle over $\Sigma$.

\end{Prop}

\begin{proof}

For each $\lambda \in U_i \cap U_j$, recall that we have

\begin{equation}
s_{\psi_j}(\lambda) = s_{\psi_i}(\lambda)\exp i \frac{\pi}{2}m_{ij}.
\end{equation}

Therefore, the above expression defines a global section if and only if 

\begin{equation}
    \frac{\Psi_j - \Psi_i}{h_k} + \frac{\pi}{2}m_{ij} + \beta_j - \beta_i = 0 \mod 2\pi
\end{equation}

The quantization condition \eqref{BSCCOHOMOLOGY} implies that

\begin{equation}
    \frac{\Psi_{j} - \Psi_{i}}{h_k} - \frac{\pi}{2} = \pi(N_k - m_k)
\end{equation}

for $j = 1$, $i = 2$ and $j = 2$, $i = 3$. Using this together with $\Psi_1 = \Psi_4$ and $\Psi_3 = \Psi_4$ on the intersection of their domains, we easily verify the values of the $\beta_j$, are determined except for a $\pm$ sign ambiguity and this is removed by requiring $\beta_2 = 0$. 
\end{proof}

\subsubsection{Conclusion of the proof of proposition 2.3}

Let$\rho_0$ be the positive invariant density on $\Sigma$ in proposition \ref{NECESSARYCONDITIONSQUANT} Define the symbol $\sigma_0 \in S^0(\Sigma)$ by

\begin{equation}
    \sigma_0 = \left( \sum_j \exp i\frac{\Psi_j(\lambda)}{h_k} \chi_j(\lambda) s_j(\lambda) \right)\sqrt{\rho_0}
\end{equation}

We inductively find a sequence of smooth functions $r_j(\lambda)$ on $\Sigma$ and complex numbers $E_j$ so that for each $n \geq 0$,

\begin{equation}
    \left(H_{h_k,c} - (1 + h^2E_1 + \cdots + h^{n+1}E_n)\right)\left( \mathscr{Q}(\sigma_0) + h\mathscr{Q}(r_1\sigma_0) + \cdots + h^n\mathscr{Q}(r_n\sigma_0) \right) \in \mathscr{O}^{n+2}
\end{equation}

With $r_0 = E_0 = 1$, the $n = 0$ case follows from formula \eqref{SUBPRINICPALTERM} and $\mathscr{L}_{X_f}\sigma_0 = 0$. Supposing it holds for $n \geq 0$, let 
$$U_n = \mathscr{Q}\left(\left(1 + \sum_{j=1}^n r_j \right)\sigma_0\right) \in \mathscr{O}^{n+2}$$ 

$$\mathcal{E}_n = 1 + \sum_{j=1}^{n}h^{j+1}E_j$$

Then with $E_{n+1}$ and $r_{n+1}$ to be determined, the function

\begin{equation}
    \left(H_{h_k,c} - \mathcal{E}_n - h^{n+2}E_{n+1}\right)\left(U_n +  h^{n+1}\mathscr{Q}(r_{n+1}\sigma_0) \right)
\end{equation}

is in $\mathscr{O}^{n+2}$ and its principal symbol is the same as the principal symbol of 

$$U_n + h^{n+2}E_{n+1}\mathscr{Q}(\sigma_0) + h^{n+1}(H_{h_k,c} - 1)\mathscr{Q}(r_{n+1}\sigma_0)$$

which vanishes if and only if 

\begin{equation}
    \frac{2}{i}dr_{n+1}(X_f) + E_{n+1} + u_n = 0
\end{equation}

where $h^n u_n \sigma_0 = \sigma(U_n)$. If $E_{n+1} = - \int_{\Sigma} u_n \, \rho_0$, then proposition \ref{NECESSARYCONDITIONSQUANT} implies that there is a smooth $r_{n+1}$ which solves this equation. Now letting $r \sim 1 + \sum_{j=1}^\infty r_j h^n \sigma_0$, $v_{h_k} = \mathscr{Q}(r\sigma_0)$ satisfies $(H_{h_k,c} - E(h))v_{h_k} \in \mathscr{O}^\infty(I,\Sigma,h_k)$. Finally, to verify the pointwise asymptotics, we write $K_j = i_{\psi_j}^*\chi_j$ and observe that

\begin{equation} \label{LEADINGQUANT}
    \begin{aligned}
   \mathscr{Q}(\sigma_0) = \sum_{j=1,3} K_j(\phi)a_j(\phi)e^{i\left(\frac{\psi_j}{h_k} + \beta_j \right)} + \sum_{j=2,4} (2\pi h_k)^{-\frac{1}{2}}\int K_j(\tau)a_j(\tau)e^{i\left( \frac{\psi_j}{h_k} + \beta_j \right)} \, d\tau  \end{aligned}
\end{equation}

\begin{equation}
    a_j(\phi) = \frac{1}{\sqrt{2\pi}} \frac{1}{(1 - V(\phi))^\frac{1}{4}} = \frac{1}{\sqrt{2\pi}} \frac{\sqrt{\sin \phi}}{(\sin^2\phi - c^2)^\frac{1}{4}} \quad j = 1,3
\end{equation}

\begin{equation}
    a_j(\tau) = \frac{1}{\sqrt{\pi}} \frac{1}{\abs{V'(G_j'(\tau))}^{\frac{1}{2}}} \quad j = 2,4
\end{equation}

Notice that $|G_j''(\tau)|^{-\frac{1}{2}} = \frac{1}{\sqrt{2}}\frac{|V'(G_j'(\tau))|^\frac{1}{2}}{\abs{\tau}^\frac{1}{2}}$, so if we apply stationary phase to $j=2,4$ terms at some fixed $\phi \in (\phi_-,\phi_+)$, the amplitudes match those in the $j = 1,3$ terms. Proposition \ref{PHASEPROP} implies that the phases match as well, so using the fact that the $K_j$ are a partition of unity when lifted to $\Sigma$, we get \eqref{ALLOWEDPTWISE}. The statement \eqref{TURNINGPTWISE} is obvious since the $j = 4$ term does not have any critical points away from the projection of $U_2$ and the $j=1,3$ are supported away from $(\phi_+, 0)$. Now take the real part of $v_{h_k}$. It satisfies the equation \eqref{QMEQN} with $E(h)$ replaced by its real part. The principal symbol of $\overline{v}_{h_k}$ is 

$$\sigma(\overline{v}_{h_k})(\phi,\tau) = \overline{\sigma(v_{h_k})}(\phi,-\tau) = \sigma_0$$

so $||\text{Re}(v_{h_k})||_{L^2(I)} = 1 + O(h)$. It follows that $L^2$ normalizing $\text{Re}\, v_{h_k}$ only multiplies the lower order terms in the full symbol by a constant. And therefore the expression for the leading part of $||\text{Re}\,v_{h_k}||^{-1}_{L^2(I)}\text{Re} \, v_{h_k}$ is the same as \eqref{ALLOWEDPTWISE}.

\subsection{Comparison of the quasi-mode to the mode}

Here we show that the quasi-mode $v_{h_k}$ of proposition \ref{QMPROP} is locally uniformly close to the true mode $u_{h_k} = UP^{m_k}_{N_k}$.

\begin{Prop} \label{L2CLOSE}
Let $v_{h_k} \in \mathscr{O}^0(I,\Sigma,h_k)$ be as in proposition \ref{QMPROP} and $u_{h_k}$ be the $L^2$ normalized, real-valued function satisfying $ H_{h_k,c}u_{h_k} = u_{h_k}$. Let $\Pi$ denote orthogonal projection onto $\ker H_{h_k,c} - 1$. Then

\begin{equation}
    ||v_{h_k} - \Pi v_{h_k}||_{L^2(I)} = O(h_k^\infty)
\end{equation}

\end{Prop}

\begin{proof}
From the spectral gap \eqref{SPECGAP}, it follows that the lower bound

$$||(H_{h_k,c}- 1)u||_{L^2(I)} \geq \delta h_k||u||_{L^2(I)}$$

holds for $u \in (\ker H_{h_k,c} - 1)^\perp$. The estimate 

$$||(H_{h_k,c} - E(h))v_{h_k}||_{L^2(I)} = O(h_k^\infty)$$

implies that there is an eigenvalue of $H_{h_k}$ in an $O(h_k^N)$ neighborhood of $E(h)$ for all large $N$. Since $E(h) = 1 + O(h^2)$ and the eigenvalues of $H_{h_k,c}$ are separated by $O(h)$ distances, this means that $E(h) = 1 + O(h_k^\infty)$. Therefore

\begin{equation}
    ||(H_{h_k,c} - 1)(v_{h_k} - \Pi v_{h_k})||_{L^2(I)} = O(h_k^\infty)
\end{equation}

which proves the estimate in view of the lower bound above.

\end{proof}

\begin{Prop} \label{UNIFORMLYCLOSE}
For each $\delta > 0$, with $I_\delta = (\delta,\pi-\delta)$,

\begin{equation}
    ||v_{h_k} - u_{h_k}||_{L^\infty(I_\delta)} = O_\delta(h_k^\infty)
\end{equation}

\end{Prop}

\begin{proof}
Writing $\partial_\phi^2 = -h_k^{-2}(H_{h_k,c} - V)$. we have

$$||\partial^2_{\phi}(v_{h_k} - \Pi v_{h_k})||_{L^2(I_\delta)} \leq h_k^{-2}\left(||(H_{h_k,c}(v_{h_k} - \Pi v_{h_k})||_{L^2(I_\delta)} + ||V(v_{h_k}- \Pi v_{h_k})||_{L^2(I_\delta)} \right)$$

From proposition \ref{L2CLOSE}, $||H_{h_k,c}(v_{h_k} - \Pi v_{h_k})||_{L^2(I)} = O(h_k^\infty)$ and since $V$ is bounded on $I_\delta$ depending on $\delta$, the right hand side is $O_\delta(h_k^\infty)$. Applying the Sobolev estimate

$$||f||_{L^\infty} \leq C ||f'||_{L^2}$$

twice on the interval $I_\delta$ together with the above inequality yields 

\begin{equation}
    ||v_{h_k} - \Pi v_{h_k}||_{L^\infty(I_\delta)} = O_\delta(h_k^\infty)
\end{equation}

Similarly, we  have

\begin{equation}
||u_{h_k}||_{L^\infty(I_\delta)} \leq C ||u_{h_k}''||_{L^2(I_\delta)} = h_k^{-2}||\left( H_{h_k,c} - V\right)u_{h_k}||_{L^2(I_\delta)} = O_\delta(h_k^{-2})   
\end{equation}

so

\begin{equation}
    ||v_{h_k} - u_{h_k}||_{L^\infty(I_\delta)} \leq ||v_{h_k} - \Pi v_{h_k}||_{L^\infty(I_\delta)} + ||(\zeta(h_k) - 1)u_{h_k}||_{L^\infty(I_\delta)} = O_\delta(h^\infty_k)
\end{equation}

Where we have written $\Pi v_{h_k} = \zeta(h_k)u_{h_k}$ and $\zeta(h_k) = 1 + O(h_k^\infty)$ since $v_{h_k}$ is real valued and positive in a neighborhood of $\phi = 0$.

\end{proof}

\section{Airy expansion of $v_{h_k}$ at the turning points}

The goal of this section is to prove the following Airy expansion for $v_{h_k}$ in a neighborhood of the turning points $\phi_\pm$. 

\begin{Prop} \label{LEGAIRY}
Let $v_{h_k}$ be the quasi-mode in proposition \ref{QMPROP}. There exists $\epsilon > 0$ such that for $\phi_+ - \epsilon < \phi < \phi_+$, $v_{h_k}(\phi)$ has the full asymptotic expansion 

\begin{equation}
    v_{h_k}(\phi) \sim Ai\left( -h^{-\frac{2}{3}} \rho(\phi) \right)h^{-\frac{1}{6}}_k \sum_{n = 0}^\infty u_{0,n}(\phi)h^{n} + Ai'\left(-h^{\frac{2}{3}}\rho(\phi)\right)h^{\frac{1}{6}}\sum_{n=0}^\infty u_{1,n}(\phi)h^n 
\end{equation}

The leading part of the expansion is 

\begin{equation}
    v_{h_k}(\phi) \sim \sqrt{\sin\phi} \, h_k^{-\frac{1}{6}} \left(\frac{4\rho(\phi)}{\sin^2\phi - c^2}\right)^\frac{1}{4}Ai\left(-h_k^{-\frac{2}{3}}\rho(\phi)\right) + O(h_k^{\frac{1}{6}})
\end{equation}

Here, the argument of the Airy function is 

\begin{equation} \label{AIRYARG}
    \rho(\phi) = \left(\frac{3}{4}\int_{\gamma_\phi} \alpha\right)^{\frac{2}{3}}
\end{equation}

where $\gamma_\phi$ is the arc on $\Sigma$ passing through $(\phi_+,0)$ from $(\phi,\tau_-)$ to $(\phi,\tau_+)$.

\end{Prop}

To prove this we write 

\begin{equation} \label{LOCALTURNINGPOINT}
    v_{h_k}(\phi) = (2\pi h_k)^{-\frac{1}{2}}\int a(\tau,h) \exp i\left( \frac{\psi_4(\phi,\tau)}{h_k} + \beta_4 \right) + O(h_k^\infty)
\end{equation}

For $\phi$ in a neighborhood of the turning point $(\phi_+,0)$. The expansion is a consequence of the following proposition from H\"ormander:

\begin{Prop}[Ho1, Theorem 7.7.18] \label{AIRYEXPANSION}
Let $f(t,x)$ be a real-valued smooth function defined in a neighborhood $(0,0) \in V \subset \R^2$. Suppose that $\partial_t f(0,0) = \partial^2_t(0,0)$ and $\partial_t^3f(0,0) \neq 0$. Then there exists smooth, real-valued functions $a(x),b(x)$ and smooth compactly supported functions $u_{0,n}(x)$, $u_{1,n}(x)$ such that 

\begin{equation}
\begin{aligned}
    e^{-\frac{i}{h}b(x)}\int u(t,x) e^{\frac{i}{h}f(t,x)} \, dt \sim & Ai(h^{-\frac{2}{3}} a(x)) h^{\frac{1}{3}}\sum_{0}^\infty u_{0,n}(x)h^{n} \\
    &+ Ai'(h^{-\frac{2}{3}}a(x)) h^{\frac{2}{3}} \sum_{n = 0}^\infty u_{1,n}(x)h^n
\end{aligned}
\end{equation}

For a smooth, compactly supported amplitude $u(t,x)$ supported sufficiently close to $(0,0)$.

\end{Prop}

\subsection{Proof of Proposition \ref{LEGAIRY}}

As explained in \cite{LPDO1}, page 234 the functions $a(x)$ and $b(x)$ can be calculated by putting the phase function into the following cubic normal form

\begin{Prop}[Ho1, Theorem 7.5.13] \label{CUBICNORMALFORMPROP}
Let $f(t,x)$ be a real valued smooth defined in a neighborhood  $(0,0) \in V \subset \R^2$ such that $\partial_t f(0,0) = \partial_t^2 f(0,0) = 0$ and $\partial_t^3(0,0) \neq 0$. Then there exists a real valued smooth function $T(t,x)$ in a neighborhood of $(0,0)$ with $T(0,0) = 0$, $\partial_t T(0,0) > 0$ and smooth functions $a(x),b(x)$ such that

\begin{equation}
    f(t,x) = \frac{T^3(t,x)}{3} + a(x)T(t,x) + b(x)
\end{equation}

\end{Prop}

We apply this theorem to the phase

$$\psi_4(\phi,\tau) = \phi\tau - G_4(\tau)$$

It has a degenerate critical point at the turning point$(\phi_+,0)$. Indeed, by differentiating the Eikonal equation,

\begin{equation}
    \tau^2 + \frac{c^2}{\sin^2 G'_2(\tau)} = 1
\end{equation}

We see that $\partial_\tau^2\psi_4(\phi_+,0) = -G_4''(0) = 0$ and $\partial_\tau^3\psi_4(\phi_+,0) = -G_4'''(0) = \frac{4c}{\sqrt{1-c^2}} \neq 0$. The functions $a(x)$ and $b(x)$ are calculated in the next proposition.

\begin{Prop}
    There exists a smooth function $T(\phi,\tau)$ in a neighborhood of $(\phi_+,0)$ as in proposition \ref{CUBICNORMALFORMPROP} such that
    
    \begin{equation} \label{PHASENORMALFORM}
        \psi_2(\phi,\tau) = \frac{T^3(\phi,\tau)}{3} + a(\phi)T(\phi,\tau) + b(\phi)    
    \end{equation}

If $\phi_+ - \epsilon < \phi < \phi_+$, then

\begin{equation}
    a(\phi) = - \left( \frac{3}{4}\int_{\gamma_\phi} \alpha \right)^{2/3}
\end{equation}

\begin{equation}
    b(\phi) = \beta_4
\end{equation}
    
Where $\gamma_\phi$ is the arc on $\Sigma$ defined in proposition \ref{LEGAIRY}.
    
\end{Prop}

\begin{proof}

Existence follows from proposition \ref{CUBICNORMALFORMPROP}. Put $\rho(\phi) = -a(\phi)$. Take the $\tau$-derivative of \eqref{PHASENORMALFORM} and observe that $\partial_\tau \psi_4(\phi,\tau) = \phi - G_4'(\tau) = 0$ if and only if $T^2(\phi,\tau) = \rho(\phi)$. For a fixed $\phi \in (\phi_+ - \epsilon, \phi_+)$ let $\tau_\pm$ be the two $\tau$-critical points of $\psi_4$, the $\tau$-coordinates of the two points $(\phi,\tau_\pm) \in \Sigma$ lying over $\phi$,

\begin{equation}
    \tau_\pm(\phi) = \pm \sqrt{1 - \frac{c^2}{\sin^2\phi}}
\end{equation}

Since $T^2(\phi,\tau_\pm(\phi)) = \rho(\phi)$, we can write $T(\phi,\tau_{+}(\phi)) = -\sqrt{\rho(\phi)}$ and $T(\phi,\tau_-(\phi)) = \sqrt{\rho(\phi)}$. These imply that $\psi_2(\phi,\tau_\pm) = \mp\frac{\rho(\phi)^{3/2}}{3} - \pm \rho^{3/2}(\phi) + b(\phi)$ which means

\begin{equation}
\frac{4}{3}\rho^{3/2}(\phi) = \psi_4(\phi,\tau_+) - \psi_4(\phi,\tau_-) \quad \quad \quad 2b(\phi) = \psi_4(\phi,\tau_+) + \psi_4(\phi,\tau_-)
\end{equation}

The formulas then follow since $\psi_4$ is odd in $\tau$ and $\Psi_4(\tau) = \psi_4(G_4'(\tau),\tau)$ is a primitive for $\alpha|_{U_4}$.

\end{proof}

Now let $\chi(\phi)$ be a bump function equal to one on $(\phi_+ - \frac{\epsilon}{2}, \phi_+ + \frac{\epsilon}{2})$ and supported in $(\phi_+ - \epsilon , \phi_+ + \epsilon)$. The amplitude $\chi(\phi)a_4(\tau,h)$ appearing in \eqref{LOCALTURNINGPOINT} will then have no critical points outside of a $\tau$ neighborhood $B_r(0)$ of $\tau = 0$, $r = o(\epsilon)$. Split up the integral by inserting a $\tau$ cutoff $\eta(\tau)$,
$\chi(\phi)a_4(\tau,h) = \chi(\phi)\eta(\tau)a(\tau,h) + \chi(\phi)(1 - \eta(\tau)a_4(\tau,h)$ supported on $B_r(0)$, equal to $1$ on $B_{r/2}(0)$. If $\epsilon$ is small enough, the first term is supported close enough to $\tau = 0$ to apply proposition \ref{AIRYEXPANSION}, and the second term is $O(h_k^\infty)$. Finally, we calculate the leading order amplitude $u_{0,0}(\phi)$ appearing in the expansion. The leading term is 

\begin{equation}
    v_{h_k}(\phi) \sim (-1)^{m_k + N_k} (2\pi)^{\frac{1}{2}}h_k^{-\frac{1}{6}} u_{0,0}(\phi) Ai(-h_k^{-\frac{2}{3}}\rho(\phi))
\end{equation}

Using the standard expansion of the Airy function for large $t > 0$ (see \cite{LPDO1} pg. 215)

\begin{equation}
    Ai(-t) \sim \frac{1}{\sqrt{\pi}t^{1/4}} \cos(\frac{2}{3}t^{\frac{3}{2}} - \frac{\pi}{4})
\end{equation}

we see that when $h^{-\frac{2}{3}}\rho(\phi) >> 0$, 

\begin{equation}
    v_{h_k}(\phi) \sim u_{0,0}(\phi)\left(\pi^{-\frac{1}{2}}  \rho(\phi)^{-\frac{1}{4}}\sin \left( h_k^{-1}\int_{\gamma_\phi}\alpha + \frac{\pi}{4} \right)\right)
\end{equation}

But this must match the leading term in proposition \ref{QMPROP} which forces 

\begin{equation}
    u_{0,0}(\phi) = \sqrt{\sin\phi}\left(\frac{4\rho(\phi)}{\sin^2\phi - c^2}\right)^\frac{1}{4}
\end{equation}

\section{Geometry of ladder sequences of spherical harmonics}

Recall that the generator of rotations, $D_\theta = -i\partial_\theta$ commutes with the Laplacian on $S^2$. The Clairaut integral, $p_\theta(x,\xi) = \langle \xi , \partial_\theta \rangle_x$ is the symbol of the $D_\theta$ so $\{p_\theta , q \} = 0$ where $q(x,\xi)  = \abs{\xi}_x$. Together they generate a homogeneous Hamiltonian torus action, $\Phi(t,s)$ on $T^*S^2$,

\begin{equation}
    \Phi(t,s,x,\xi) = \exp sX_{p_\theta} \circ \exp tX_{q}(x,\xi)
\end{equation}

which acts transitively on the level sets of the moment map,

\begin{equation}
    \mu : T^*S^2 \to \Gamma \subset \R^2 \quad \quad \mu(x,\xi) = (q(x,\xi),p_\theta(x,\xi))
\end{equation}

whose image is the cone $\Gamma = \{(x,y) ~|~ x \geq 0, \abs{y} \leq x\}$. Since $\mu$ is homogeneous, we need only consider level sets with $q = 1$. For $c \in [-1,1]$ set $T_c = \mu^{-1}(1,c) = S^*S^2 \,  \cap \, \{ p_\theta = c \}$. For $c \in (-1,1)$, $T_c$ is Lagrangian torus inside of $S^*S^2$. When $c = \pm 1$, $T_c$ is just the lift of the standard equator $\gamma_e$ to $S^*S^2$. In terms of the polar coordinates $(\phi,\tau,\theta,\theta)$ on $T^*S^2$, 

\begin{equation}
    T_c = \{(\phi,\tau,\theta,\eta) ~|~ \tau^2 + \frac{c^2}{\sin^2\phi} = 1 ; \eta = c \}
\end{equation}

Therefore its projection to $S^2$ is $\pi(T_c) = \{ (\phi,\theta) ~|~ \sin \phi \geq c \}$. The projection is an annular band consisting of all geodesics which make the fixed angle $\arccos c$ with the equator. The energy curve associated with the associated Legendre functions is just the $(\phi,\tau)$ cross-section of $T_c$. For $x$ in the interior of $\pi(T_c)$, let $\gamma_x$ be the geodesic arc connecting the two points lying above $x$ in $T_c$, from $(x,\xi_-)$ to $(x,\xi_+)$, where the sign corresponds to the sign of $\tau$. It is clear that $\int_{\gamma_x} d\theta = 0$ since there is no change in the $\theta$ coordinate across the arc. But the canonical 1-form is $\alpha = \tau d\phi + \eta d\theta$. Since $\eta = c$ is constant on $T_c$ the second term contributes nothing to the integral over $\gamma_x$, and the first term is clearly equal to the integral \eqref{AIRYARG} in the Legendre function expansion. The density 

\begin{equation}
     d\mu_c = \frac{ \abs{d\phi} \otimes \abs{d\theta}}{(2\pi)^2\abs{\tau}}
\end{equation}

is invariant under the joint flow on $T_c$ and 

$$\pi_*d\mu_c = \frac{1}{(2\pi)^2} \frac{\sqrt{\sin \phi} ~\abs{d\phi}\otimes \abs{d\theta}}{\sin^2\phi - c^2}$$

which verifies the formula \eqref{HALFDENSITYINTERP}. The reason for the Airy bump at the caustic latitude circles is the presence of a fold singularity for the projection $\pi|_{T_c} : T_c \to S^2$. 
Recall that a smooth map $f : X^n \to Y^n$ between $n$-dimensional manifolds is said to have a fold singularity with fold locus $S$ if there exists a codimension one submanifold $S \subset X$ such that 

\begin{enumerate}
    \item $S$ is equal to the set of critical points of the map $f$, i.e. $S = \{ x \in X ~|~ df_x ~\text{is not surjective}~ \}$
    
    \item For each $s \in S$, the kernel of $df_s$ is transverse to the tangent space $T_s S$.
\end{enumerate}

\begin{Prop}
The projection $\pi|_{T_c} \to S^2$ is a folding map with fold locus $S = S_+ \cup S_-$,

$$S_{\pm} = \{(\phi_{\pm}, \theta, 0 , c) ~|~ \theta \in [0,2\pi) \}$$

where $\phi_\pm$ are the two solutions of $\sin \phi = c$. The images of $S_\pm$ are the latitude circles which bound $\pi(T_c)$.

\end{Prop}

\begin{proof}

Writing $(\rho,\eta)$ as dual coordinates to $(\phi,\theta)$, $T_c$ is cut out by the equations $\eta = c$ and $\rho^2 + \frac{c^2}{\sin^2\phi} = 1$. differentiating the second equation gives

$$\rho d\rho - c^2 \frac{\cos\phi}{\sin^3\phi}d\phi = 0 $$

so writing $x = (\phi,\theta)$, $\xi = (\rho,\eta)$,

$$T_{(x,\xi)}T_c = \{ \alpha\partial_\phi + \beta\partial_\theta + \gamma\partial_\rho ~|~ \rho\gamma = c^2\frac{\cos \phi}{\sin^3 \phi}\alpha \}$$

So for $v \in T_(x,\xi)T_c$, $d\pi v = 0$ if and only if $\alpha = \beta = 0$. But then $\rho\gamma = 0$. If $\rho = 0$, then $v = 0$, so the kernel of $d\pi$ is non-trivial only when $\rho = 0$, and this means that $(x,\xi) \in S$. At such points, the kernel of $d\pi$ is the span of $\partial_\rho$, which is transverse to $T_{(x,\xi)}S = \R{}\partial_\theta$.

\end{proof}

\section{Calculation of the weak limit of the empirical measures, proof of theorem 1.2}

In this section we determine the weak limit of the measures $\eqref{EMPIRICALMEASURE}$. Noting that $\abs{Y^m_N(x)}^2$ is constant on latitude circles, choose any $x \in \gamma_{c_0}$ and we may write

\begin{equation}
    \mu_{N,c_0}(t) = \frac{2\pi}{N + \frac{1}{2}} \,\sum_{m = -N}^N \abs{Y^m_N(x)}^2\delta_0\left(t - \frac{m}{N}\right).
\end{equation}

Thus $\int \mu_{N,c_0} = \frac{2\pi}{N + \frac{1}{2}} \Pi_{N}(x,x) = 1$. By the L\'evy continuity theorem, it suffices to show that the characteristic functions

\begin{equation}
\Phi_{N,c_0}(s) = \frac{1}{2\pi} \int e^{is t} \, d\mu_{N,c_0}(t)    
\end{equation}

converge pointwise to a limit $\Phi(s)$ which is continuous at $s = 0$. In this case, $\mu_{N,c_0}$ converges weakly to the Fourier transform of $\Phi(s)$. 

\begin{Prop}

The characteristic functions $\Phi_{N,c_0}$ converge pointwise to $(2\pi)^{-1} J_0(c_0 \, s)$ where

\begin{equation}
    J_0(s) = \frac{1}{2\pi} \int_{-\pi}^\pi \exp \left(-i s \sin \theta \right) \, d\theta
\end{equation}

is the order zero Bessel function.

\end{Prop}

\begin{proof}

Write 

\begin{equation}
    \Phi_{N,c_0} = \frac{1}{N + \frac{1}{2}} \sum_{m = -N}^N \exp\left( i \frac{sm}{N} \right) \abs{Y^m_N(x)}^2
\end{equation}

Let $r_s : S^2 \to S^2$ be rotation by $s$ in the polar angle, i.e. pullback under the flow of $\partial_\theta$. Then the right hand side can be rewritten as

\begin{equation}
    \Phi_{N,c_0}(s) = \frac{1}{N + \frac{1}{2}} \Pi_N(r_{s/N}(x),x) = \frac{1}{2\pi} P_N\left( \cos d(r_{s/N}(x),x) \right)
\end{equation}

where $d(x,y)$ is the Riemannian distance. Recall $\cos d(x,y) = x \cdot y$ so $d(r_s(x),x) = 1 - c_0^2(1 - \cos s)$ and in particular $d(r_s(x),x) \leq c_0\abs{s}$, hence

\begin{equation} \label{DISTANCEASYMPTOTICS}
\begin{aligned}
 \frac{d(r_{s/N}(x),x)^2}{2} &= 1 - \cos d(r_{s/N}(x),x) + O\left(\frac{s^4}{N^4}\right) \\
 &= c_0^2 \frac{s^2}{2N^2} + O\left(\frac{s^4}{N^4}\right)
\end{aligned}
\end{equation}

By the Mehler-Heine asymptotics we have, locally uniformly in $z$, 

\begin{equation}
    \lim_{N \to \infty} P_N(\cos {\frac{z}{N}}) = J_0(z)
\end{equation}

which, together with \eqref{DISTANCEASYMPTOTICS} implies the limit.

\end{proof}

Now the weak limit calculation is finished in light of the fact that 

\begin{equation}
    \int_{-\infty}^{\infty} e^{-ist}J_0(s) \, ds = 1|_{[-1,1]}\frac{2}{\sqrt{1 - t^2}}
\end{equation}

which is easily verified directly.

\end{document}